\documentclass[a4paper]{article}
\usepackage{graphicx}  		
\usepackage[pdfauthor={Amaro Rica da Silva},%
pdftitle={Galois Groups in the work of Mira Fernandes},%
pagebackref=true,%
pdftex]{hyperref}
\usepackage{amsmath, amssymb}
\usepackage{endnotes}
\usepackage[latin1]{inputenc}
\def\dn#1#2{	\raisebox{-#1pt}{$#2$}	}	
\def\Zeta{	\zeta\dn2		}
\newcommand{\csep}[2][1.4]{\renewcommand{\arraystretch}{#1}\setlength\arraycolsep{#2 bp}}

\def\scP{{\cal P}}
\def\scF{{\cal F}}
\def\scR{{\cal R}}
\def\scQ{{\cal Q}}
\def\scI{{\cal I}}
\newcounter{nummer} 
\newcommand{\fn}[2]{\hypertarget{back\arabic{nummer}}{\hyperlink{note#1}{\endnotemark}}%
\stepcounter{nummer}%
\endnotetext{\protect\hypertarget{note#1}{#2}}}

\begin{document}

\title{Galois Groups in the work of Mira Fernandes\\
}
\author{Amaro Rica da Silva \\
{\small Centro Multidisciplinar de {}Astrof\'{\i}sica -
CENTRA,}\\
{\small Departamento de F\'{\i}sica, Instituto Superior
T\'ecnico - IST,}\\
{\small Universidade T\'ecnica de Lisboa - UTL,}\\
{\small Avenida Rovisco Pais 1, 1049-001 Lisboa, Portugal,}}
\date{
}
\maketitle
\setlength{\parindent}{0em}


Aureliano Mira Fernandes was a student at the University of Coimbra from 1904 until 1910 when he finished his Mathematics degree. He studied Calculus with Sidónio Pais and Analysis with José Bruno de Cabedo.
In March, 1911, he completed his Ph.D. thesis  entitled "Teorias de Galois
I-Elementos da teoria dos grupos de
substituições"\footnote{"Galois Theories I - Elements of
the theory of finite substitution groups"
reedited as "Substitution Groups and Algebraic Solvability I" (1929)\cite{AdMF29a}}  under the orientation of Prof. Souto Rodrigues of the University of Coimbra. He was then invited to be Full Professor at the IST in November the same year 1911.

His Ph.D. thesis was presented at a time when group theory ideas
were emerging at the forefront of scientific research in many areas in most European countries but were not widely known. For instance the first english expositions and translations of Galois theory appear around 1891-1900 but are essentially geared towards construction methods of Galois groups (O. Bolza, J. Pierpoint, H. Voigt) and up to 1908 no course on Galois theory was taught at Cambridge or Oxford.

As was kindly noted to me by Prof. Paulo Almeida, there were a few Portuguese mathematicians in the late 1900 that studied algebraic equations, such as Prof. Luiz Woodhouse at the Academia Politécnica do Porto, who included the subject in his course Higher Algebra and Analytical Geometry, and the Jornal de Sciencias Mathematicas e Astronomicas da Universidade de Coimbra published works on algebraic equations by Martins da Silva (1882) and Whoodhouse (in 1885) where reference to results by Galois can be found, but these works do not address Galois methods or theory. 

In his 1911 Ph.D. thesis  Mira Fernandes focuses on the results around which the group structure behind the solvability theorems of Galois theory reside, and deals mainly with finite group theory definitions and results with permutation group realizations. 
Although he refers to algebraic equation root finding and Galois theory he does not deal with these methods in his thesis, which he leaves until a later publication in 1931, in part II of a work entitled "Grupos de Substituições e Resolubilidade Algébrica" published by the Instituto Superior de Comércio de Lisboa. 
The latter publication is definitely of a pedagogical nature by the time it is published, as Mira Fernandes' interests by then had evolved to the applications of Lie Groups to general relativistic theories. 

The following table of contents from these works illustrates the subject matter dealt with in each publication.

\begin{minipage}[r]{\hsize}\footnotesize
\subsubsection*{Galois Theories I - Elements of the theory of finite substitution groups (1911)\cite{AdMF08}}
\begin{enumerate}
\item[] Introduction: Algebraic Solvability
\item[{I}- ]Finite Groups: Transitivity and primitivity
\item[ {II}- ]Isomorphism and Group composition: Jordan-Hölder and Sylow theorems, Solvable groups.
\item[ {III}- ]Abelian Groups
\item[ {IV}- ]Metacyclic group, General Linear group and the Modular group.
\item[ {V}- ]The structure of the total group and the Alternating group - Possible orders of the simple groups.
\item[ {VI}- ]Generalization of the concept of isomorphism: linear substitution groups of finite order.\footnotemark[2]
\item[ {VII}- ]Geometrical representation of finite groups of linear substitutions: groups of regular polyhedra.
\end{enumerate} 

\subsubsection*{Substitution Groups and Algebraic Solvability II (1931)\cite{AdMF29b}}
\begin{enumerate}
\item[ {I}- ]Algebraic Field-Irreducibility of polynomials
\item[ {II}- ]Galois Resolvent. Galois Group and properties
\item[ {III}- ]General Resolvent. Structure of the Galois group.
\item[ {IV}- ]Abelian Equations.
\item[ {V}- ]Solvability via radicals.
\end{enumerate}
\end{minipage} 
\footnotetext[2]{meaning $z\to z'=\frac{p\, z+q}{p'\,z+q'}$ with
$p\,q'-p'\,q\neq 0$ }
\setcounter{footnote}{2}

By the end of the XIX century, and definitely by  1931 the focus of Galois Theories had shifted towards Number Theory applications by the exploration of extension fields mostly by the German mathematicians following the work of Dedekind. Nowadays one can see that this path to Number Theory can be summarized with the following diagram\vspace{-12pt}
\begin{center}
\begin{tabular}{ll}
$\hskip -14pt \left(\!\!
\begin{array}{c}
\mbox{Solution of Algebraic Equations}\\
\mbox{Geometric Construction Problems}
\end{array}\!\! \right) $
$ \Longrightarrow $ 
$ \left(\!\!
\begin{array}{c}
\mbox{Theory of} \\
\mbox{Polynomials}
\end{array}\!\! \right) $
 $ \Longrightarrow $ 
$ \left(\!\!
\begin{array}{c}
 \mbox{Commutative} \\
\mbox{Field Theory} 
\end{array}\!\! \right) $
\end{tabular}
\end{center}

In this presentation we intend to show how the group theory concepts discovered by Galois are connected to previous work that took centuries to evolve, and how fast they changed the landscape of theoretical and applied mathematics and physics since its publication. 
It was the concept of radicals and their use in the solution of algebraic problems that ultimately led to the concept of fields and Galois group theory. Afterwards a revolution took place as his group theory concepts took hold with applications in great many areas in mathematics and physics.
\nocite{AdMF08}\nocite{JD78}\nocite{JT01}\nocite{MA08}\nocite{MK72}\nocite{NB74}

The road that leads to Galois theory is made of contributions by mathematicians that were trying to find ways to express the general roots of algebraic equations using rational expressions and radicals involving the constant coefficients of these equations, as was made by del Ferro, Cardano and Ferrari in the XVI$\rm  {}^{th} $ century for equations of order up to 4. 
\fn{1}{\textbf{\large Algebraic equations until the XVI$\rm  {}^{th} $ century}
\setlength{\parindent}{0in}\medskip\newline
We can trace back the (numerical) methods for the solution of quadratic and bi-quadratic equations almost $ 4\ 000 $ years as it clear that Babylonians already knew how perform square-root operations since ca. $ 1900 $ BD. The following cuneiform tablet shows how to compute $ \sqrt{2} $ and $ \frac{1}{\sqrt{2}} $ solving what is basically a Pythagorean problem numerically. Babylonian mathematics is a set of numerical recipes for solving day-to-day problems, and even though there are many examples of training exercises for apprentices, there was never an attempt to formalize the theory.
\medskip\newline
\begin{minipage}{0.5\linewidth} 
\includegraphics[scale=.5]{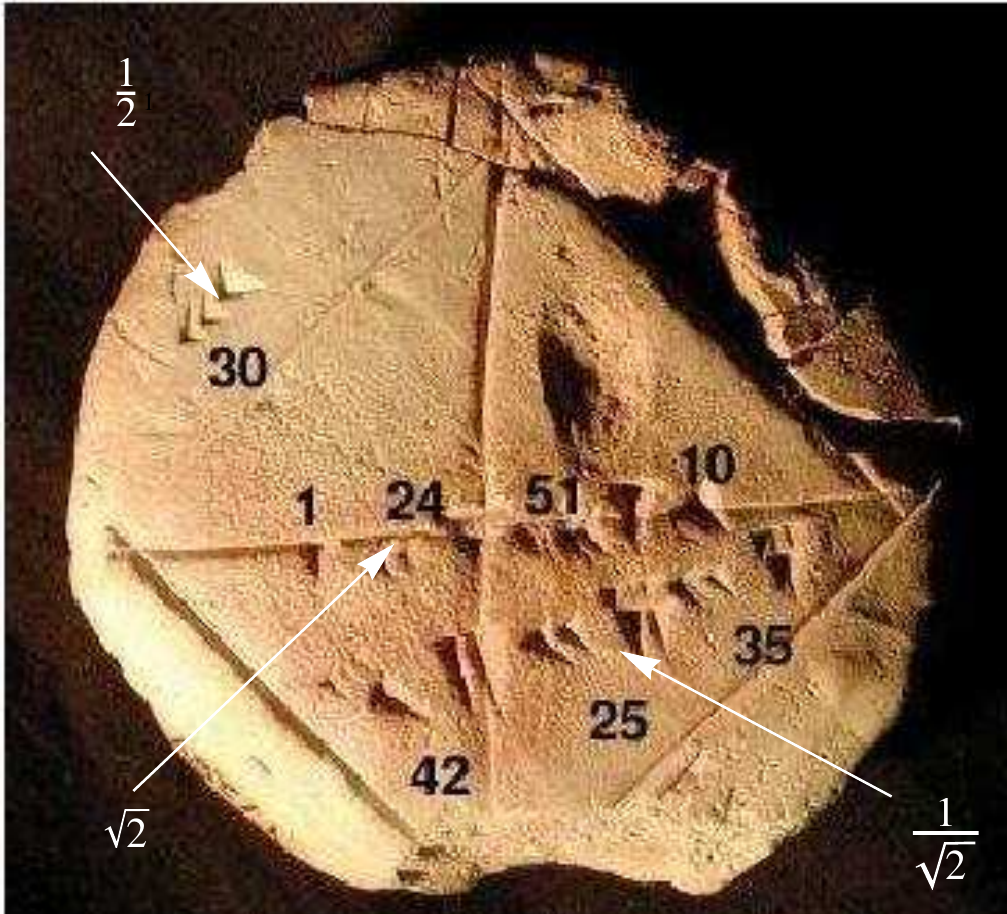}
\end{minipage}
\begin{minipage}{0.5\linewidth}
\begin{equation*}\renewcommand{\arraystretch}{2}
\left\lbrace\begin{array}{rl}
0+\frac{30}{60}&=\frac{1}{2}\\
0+\frac{42}{60}+\frac{25}{60^{2}}+\frac{35}{60^{3}}&=0.707106\underline{48148}\\
\dfrac{1}{\sqrt{2}}&=0.70710678118\\
1+\frac{24}{60}+\frac{51}{60^2}+\frac{10}{60^3} & =  1.41421\underline{296296} \\
\sqrt{2} & =  1.41421{356237}\\
\end{array}\right.
\end{equation*}
\end{minipage}
\medskip\newline
The next relevant abstraction level is introduced by the Greek geometry.
Greeks use Babylonian formulas through geometrical constructs only, and refrain from using algebraic formulas until $ 100 $ AD. (Heron and Diophantus). They did however develop methods of infinite approximation to square-roots by the IV$\rm  {}^{th} $ century AD (Théon of Alexandria).
\newline
Greek dedication to geometrical methods is probably a consequence of the fact that, in the "Elements", Euclid restricts himself to using straight edge and compass methods only, as suggested by Plato. 
\newline
Greek mathematicians also used geometrical methods to represent incommensurate ratios, which they did not consider numbers, and knew already that an irreducible algebraic equation of third degree  over the field of rationals, such as the duplication of the cube ($x^3=2$) or the trisection of an angle ($ x^{3}-3\,x-2\,b=0 $), cannot be solved with straight edge and compass only. In fact we know today that few algebraic equations possess roots that may be  found in this way.
\newline
Leonardo de Pisa is responsible for the introduction of Arabian mathematical methods in the western world in the XIII$\rm {}^{th} $ century, and the solution of quadratic equations is then perfected through formulas using radicals.
}
\fn{2}{\textbf{\large Scipione del Ferro (Bologna, 1462-1526)}
\setlength{\parindent}{0in}
\medskip\newline
In  $ 1512 $ del Ferro solves the general cubic equation.A general cubic equation $u^3+a_2\, u^2+a_1\,u+a_0=0$ can always be reduced to a form ${x^3+a\, x=b}$ by the substitution $u = x-\frac{1}{3}\,a_2$ with:
\begin{equation*}
a=a_1-\frac{{a_2}^2}{3}\qquad;\qquad b=-a_0+\frac{a_1\, a_2}{3}-\frac{2\, {a_2}^3}{27} 
\end{equation*} 
Making a variable substitution $x= y-z$ on $x^3+a\, x=b$ yields
\begin{equation*}
y^3-z^3+(y-z)\, (a-3\, y\, z) =b  
\end{equation*} 
Under the conditions for 
\begin{equation*}\renewcommand\arraystretch{1.5}
(y-z) (a-3\, y\, z)=0 \qquad\Longleftrightarrow\qquad
\left\lbrace
\begin{array}{rc}
 y\, z&=\frac{a}{3} \\
 y    &=z\\
\end{array}\right.
\end{equation*} 
del Ferro solves the two-variable system for $ y^{3} $ and $ z^{3} $
\begin{equation*}\renewcommand\arraystretch{1.5}
\left\lbrace
\begin{array}{rl}
 y^3-z^3 & =  b \\
 y^3\, z^3 & =\left(\frac{a}{3}\right)^3\\
\end{array}\right.
\end{equation*} 
thus obtaining
\begin{equation*}
 y^3=\frac{b}{2}+{\textstyle\sqrt{\left(\frac{a}{3}\right)^3+\left(\frac{b}{2}\right)^2}} \qquad ;\qquad
  z^3=-\frac{b}{2}+\textstyle\sqrt{\left(\frac{a}{3}\right)^3+\left(\frac{b}{2}\right)^2}
\end{equation*} 
\medskip\newline
del Ferro never publicizes his findings, sharing it only with one of his students, which latter will divulge it in the form of a sonnet.
}
\fn{3}{\textbf{\large Jeronimo Cardano (Pavia,1501-1543)}
\setlength{\parindent}{0in}
\medskip\newline
\begin{minipage}{0.2\linewidth}
 \vspace{-5pt}\includegraphics[scale=.38,keepaspectratio=true]{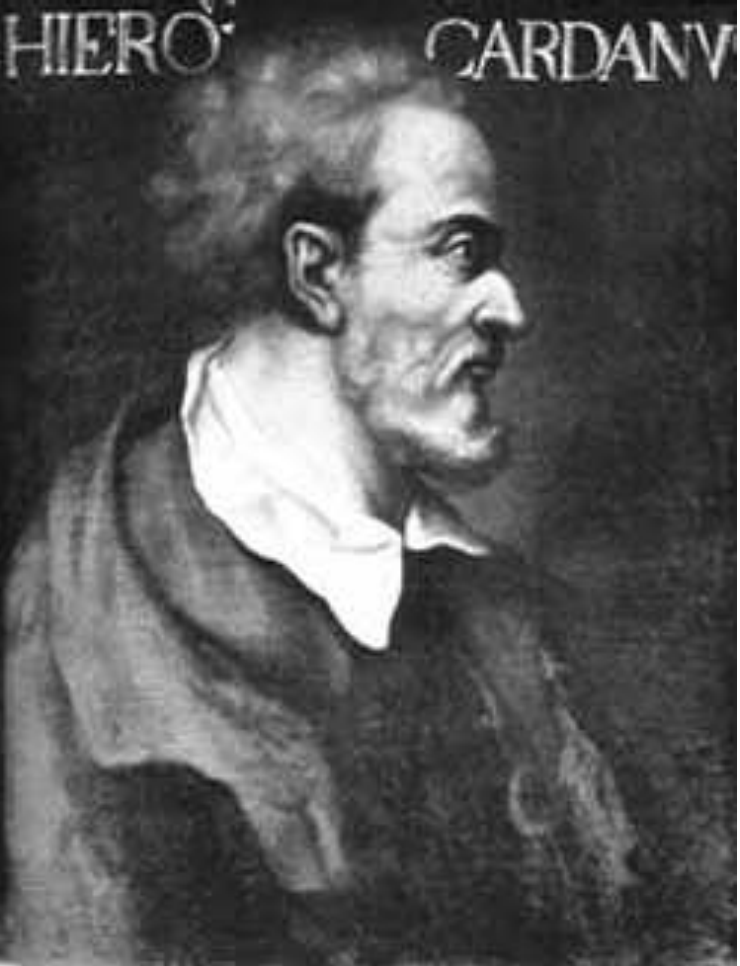}
\end{minipage}\hspace{.05\linewidth}
\begin{minipage}{0.75\linewidth}
 In $ 1545 $ Cardano publishes in his "Ars Magna" the del Ferro-Cardano formula for the cubic $x^3+a\, x=b.$
\begin{equation*}
x= \left(\frac{b}{2}+\sqrt{\frac{a^3}{27}+\frac{b^2}{4}}\right)^{1/3}- \left(-\frac{b}{2}+\sqrt{\frac{a^3}{27}+\frac{b^2}{4}}\right)^{1/3}
\end{equation*} 
There was at the time great resistance to this formula, not only because it used negative numbers, but also because real roots had to
be obtained by summing what we now know as complex numbers. Cardano never uses the complex solutions that are implicit in the radicals of $ y^{3} $ and $ z^{3} $ when considering the 3-roots of unity $\Zeta{{}_k}^{\!\!(3)}=e^{2\,\pi\,i\,\frac{k}{3}}=(-1)^{\frac{\scriptstyle 2\,k}{\scriptstyle 3}}$. 
\end{minipage}
In fact
\begin{equation*}\renewcommand\arraystretch{2.0}
 y^3=\frac{b}{2}+\sqrt{\frac{a^3}{27}+\frac{b^2}{4}} \qquad\Longrightarrow\qquad       
\left\lbrace
\begin{array}{lcrl}
 y_{3}&=&&\left(\frac{b}{2}+\sqrt{\frac{a^3}{27}+\frac{b^2}{4}}\right)^{1/3} \\
 y_{2}&=&-(-1)^{1/3} &\left(\frac{b}{2}+\sqrt{\frac{a^3}{27}+\frac{b^2}{4}}\right)^{1/3} \\
 y_{1}&=&(-1)^{2/3} &\left(\frac{b}{2}+\sqrt{\frac{a^3}{27}+\frac{b^2}{4}}\right)^{1/3}\\
\end{array}\right.
\end{equation*}
\medskip\newline
\begin{equation*}\renewcommand\arraystretch{2.0}
z^3=-\frac{b}{2}+\sqrt{\frac{a^3}{27}+\frac{b^2}{4}} \qquad\Longrightarrow\qquad  
\left\lbrace\begin{array}{lcrl}
 z_{3}&=&& \left(-\frac{b}{2}+\sqrt{\frac{a^3}{27}+\frac{b^2}{4}}\right)^{1/3} \\
 z_{2}&=&-(-1)^{1/3} &\left(-\frac{b}{2}+\sqrt{\frac{a^3}{27}+\frac{b^2}{4}}\right)^{1/3} \\
 z_{1}&=& (-1)^{2/3} &\left(-\frac{b}{2}+\sqrt{\frac{a^3}{27}+\frac{b^2}{4}}\right)^{1/3}
\end{array}\right.
\end{equation*}
\medskip\newline
Only in 1732 will Euler show that of the $9$ possible combinations of $x_{ij}=y_{i}-z_{j}$ only those with $y_{i}\, z_{j}=\frac{a}{3}$ are good roots and that the correct formulas are:
\medskip\newline
\begin{equation*}
\begin{cases}
 x_k=\Zeta{{}_k}^{\!\!(3)} \left(\frac{b}{2}+\sqrt{\frac{a^3}{27}+\frac{b^2}{4}}\right)^{1/3}-{\Zeta{{}_k}^{\!\!(3)}}^2\left(-\frac{b}{2}+\sqrt{\frac{a^3}{27}+\frac{b^2}{4}}\right)^{1/3}
& \mbox{if  }a\geq  0 \\\\
 x_k=\Zeta{{}_k}^{\!\!(3)} \left(\frac{b}{2}+\sqrt{\frac{a^3}{27}+\frac{b^2}{4}}\right)^{1/3}+{\Zeta{{}_k}^{\!\!(3)}}^2\left(\frac{b}{2}-\sqrt{\frac{a^3}{27}+\frac{b^2}{4}}\right)^{1/3}
& \mbox{if  }a<0
\end{cases}
\end{equation*} 
Cardano also publishes the solution to the 4th order algebraic equation, attributing it to his student Lodovico Ferrari.
Nowadays one can see that the del Ferro-Cardano solution is the first clue that the roots of (real) algebraic equations must be obtained by "extending" the field $ \mathbb{R} $ with roots of the unity and other radicals.
\medskip\newline
Cardano gets his medical doctor degree in 1526 by the University of Padua.
In 1534 he starts lecturing Mathematics in Milan, but maintains his studies in medicine, astrology and magic.
In 1570 he is arrested by the Inquisition on charges of having drawn the horoscope of Jesus Christ. He was released but barred from giving
any more lectures.
}

The next evolutionary step would be taken by Leibnitz and Tschirnhaus by the end of the  XVII$\rm {}^{th} $ century, and it became apparent that the general equations of fifth order or greater were problematic in this respect. First came the realization that methods that were used to simplify equations of degree less than four would not work for these higher order equations. 
\fn{4}{ \textbf{\large E. W. von Tschirnhaus (1651-1708)}
\setlength{\parindent}{0in}
\medskip\newline
\begin{minipage}{0.2\linewidth}
 \includegraphics[scale=.3,keepaspectratio=true]{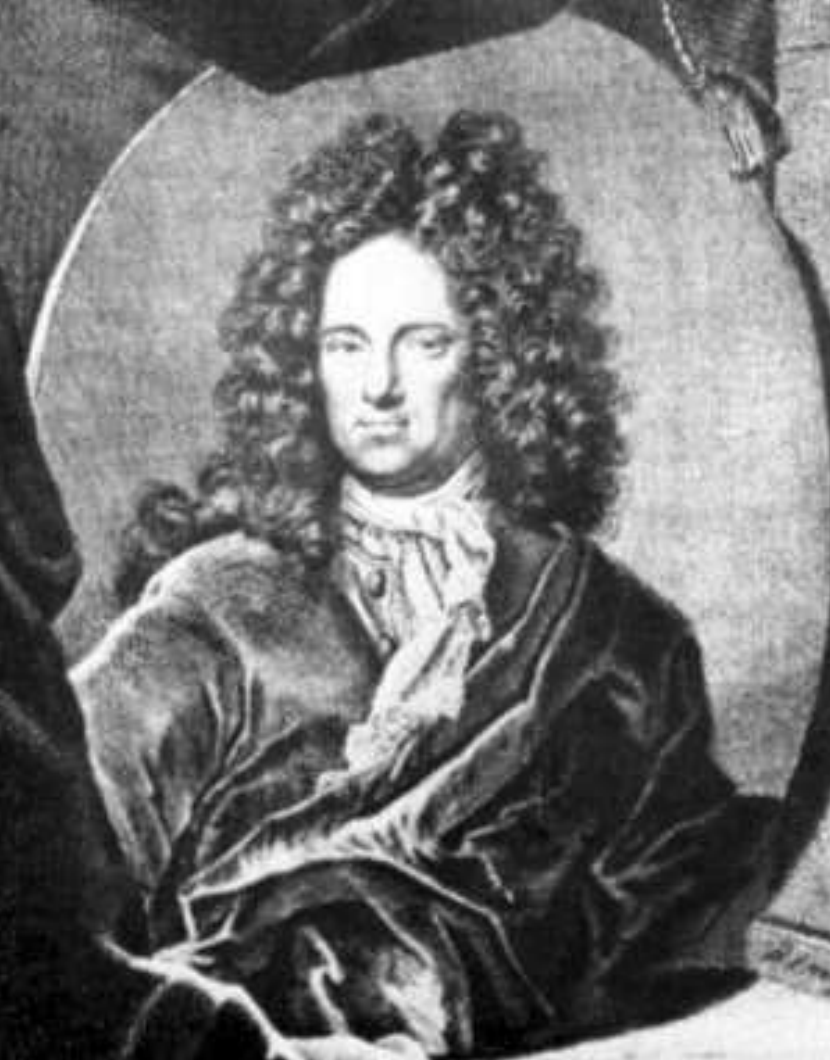}
\end{minipage}\hspace{.05\linewidth}
\begin{minipage}{0.75\linewidth}
By the end of the XVII$\rm  {}^{th} $ century Leibnitz and Tschirnhaus were among the few that still studied solutions of algebraic equations with radicals. In 1683 Tschirnhaus proposes a method to convert a general polynomial equation of degree $n$ 
\begin{equation*}
 \scP_n(x)=x^n+a_1x^{n-1}+\ldots +a_{n-1} x+a_n=0  
\end{equation*} 
into a lower degree polynomial by eliminating $x$ with an auxiliary equation of degree $n-1$ in $ x $ 
\begin{equation}
 y=x^{n-1}+b_1x^{n-2}+\ldots +b_{n-2}x+b_{n-1} 
\label{YEQ} 
\end{equation} 
\end{minipage}
\medskip\newline
Then by canceling the $ (n-1)^{2} $ coefficients of $ x^{n(n-1)},\,\dots,\,x^{n} $ in the linear combination
\begin{equation*}
\scQ_n(y)=y^n+A_1\,y^{n-1} + \ldots + A_{n-1}\,y + A_n=\left( \alpha_{n(n-1)}x^{n(n-1)}+\dots+\alpha_{1}x+\alpha_{0}\right)\scP_{n}(x)
\end{equation*} 
one can solve for the $ \alpha_{i} $ and replace these in the remaining $ n $ equations to obtain the
 $ A_{k}=A_{k}(\mathfrak{a},\mathfrak{b}) $, where $\mathfrak{a}=\left\lbrace{a_{1},\dots,a_{n}}\right\rbrace  $ and  $\mathfrak{b}=\left\lbrace{b_{1},\dots,b_{n-1}}\right\rbrace  $.
\newline
{Imposing the constraint equations }
\begin{equation}\renewcommand{\arraystretch}{1.5}
\begin{array}{lr}
 A_k(\mathfrak{a},\mathfrak{b})=0&(k=1,\,2,\dots,\,n-1) 
\end{array}
\label{CTR} 
\end{equation} 
{one could (in principle) obtain the unknown coefficients $ \mathfrak{b}=\left\{ b_{k}(\mathfrak{a})\right\} $  so that $
 A_n(\mathfrak{a},\mathfrak{b})= -c_n(\mathfrak{a}) $ and $ \scQ_n(y) $ is the binomial equation }
\begin{equation}
\scQ_n(y)=y^n-c_n(\mathfrak{a})=0 
\label{BEQ}
\end{equation} 
{With the new coefficients $b_k(\mathfrak{a})$ and the roots $y_{n,k}(\mathfrak{a})=e^{2\pi\,i\,\frac{k}{n}}\sqrt[n]{c_{n}(\mathfrak{a})}$ of this binomial equation (\ref{BEQ}) one obtains from (\ref{YEQ}) a set of lower order polynomial equations
}
\begin{equation*}
\scP_{n-1}(x,k)=x^{n-1}+b_1(\mathfrak{a})\,x^{n-2}+\ldots +b_{n-2}(\mathfrak{a})\,x+b_{n-1}(\mathfrak{a})-y_{n,k}(\mathfrak{a}) =0
\end{equation*} 
{whose roots could be determined by the same method. The goal then is to reach a point where we determine a set of monomials}
\begin{equation*}
\scP_{1}(x,k,\dots,m)=x+\beta_{1}(\mathfrak{a}) -y_{2,m}(\mathfrak{a})
\end{equation*} 
from which the roots of the original equation can be selected.
{Unfortunately this method only works for $n\leq 4$, since for $n=5 $ the constraint equations (\ref{CTR}) are of degree $24$ in the $ b_{k} $, thus being harder to solve than the original equation, and this worsens with increasing $ n $. The method also produces false roots among the genuine ones.}
}

Then in the XVIII$\rm  {}^{th} $ century Lagrange recognizes that the key to the solvability of some equations is the invariance under permutations of their arguments of certain symmetric rational functions of the roots, that he calls resolvents. In 1771 permutations were first employed by J.L. Lagrange in his "Réflexions sur la résolution algébrique des équations". By the same time Vandermonde also uses symmetric rational functions of the roots of algebraic equations (somewhat akin to Lagrange's resolvents) to find solutions by radicals of cyclotomic equations $ x^{p}-1=0 $  of degrees up to $p= 11$ (and of a particular case of an equation of degree $ 5 $ related to them).
\fn{5}{\textbf{\large Joseph-Louis de Lagrange (1736-1813)}
\setlength{\parindent}{0in}
\medskip\newline
\begin{minipage}{0.2\linewidth}
 \includegraphics[scale=.34,keepaspectratio=true]{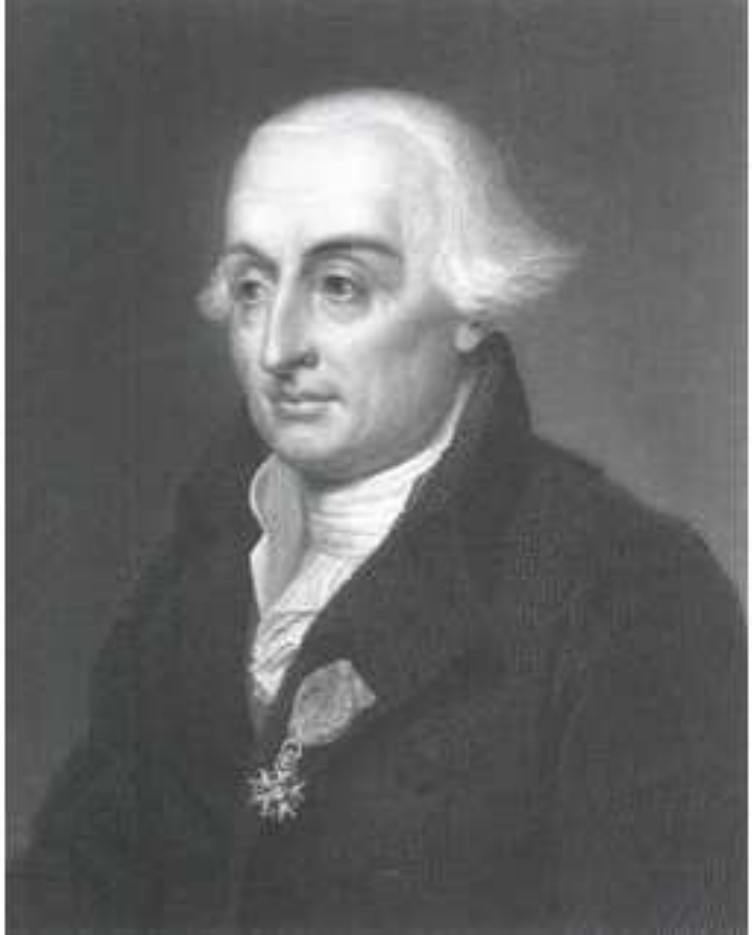}
\end{minipage}\hspace{.02\linewidth}
\begin{minipage}{0.78\linewidth}
Lagrange recognizes that what distinguishes algebraic equations of degree $n\leq 4$ for their solvability was the existence of rational symmetric functions of the roots $\beta \left(x_1,\ldots ,x_n\right)$ that under all the permutations of the $n$ roots yield a small number $m<n$ of distinct expressions $\beta _k\left(x_1,\ldots ,x_n\right)=\beta \left(x_{k_1},\ldots ,x_{k_n}\right)$. 
\newline
There is then a Resolvent Equation \begin{equation*}
                                    \scR_n(y)=y^m+b_1y^{m-1}+\ldots +b_{m-1}y+b_m=0, 
                                   \end{equation*}
with coefficients $b_k=b_k(\mathfrak{a})$ whose roots $Y_k(\mathfrak{b})=\beta _k\left(x_1,\ldots ,x_n\right)$ are precisely these invariant functions.
\end{minipage}
\medskip\newline
The Lagrange resolvents are of the form 
\begin{equation*}
 \rho _k\left(x_1,\ldots ,x_n\right)=\sum _{i=1}^n \left(\Zeta{{}_k}^{(n)}\right)^ix_i
\end{equation*} 
His invariant functions are
\begin{equation*}
 \beta _k\left(x_1,\ldots ,x_n\right)=\left(\frac{1}{n}\rho _k\left(x_1,\ldots ,x_n\right)\right)^n
\end{equation*}  
{For instance, in the cubic equation $\scP_3(x)=x^3+a x=b$, given that the 3-roots of unity obey
}
\begin{equation*}
\begin{array}{ccccc}
 \left(\zeta _1^{(3)}\right)^2=\zeta _2^{(3)} & ; & \left(\zeta _2^{(3)}\right)^2=\zeta _1^{(3)} & ; & \left(\zeta _1^{(3)}\right)^3=\left(\zeta_2^{(3)}\right)^3=\zeta _3^{(3)}=1
\end{array}
\end{equation*} 
{we obtain six different $\rho _k$ but only two distinct $\beta _k$ since the last one is identically zero.}
\begin{equation*}
 \csep[1.3]{2}
\hskip -12pt\begin{array}{|c|c|c|}\hline
k&\rho_{k}(x_{i_{1}},x_{i_{2}},\dots,x_{i_{n}})&{\rho_{k}}(x_{i_{1}},x_{i_{2}},\dots,x_{i_{n}})^{3}\\\hline
 1,2	&	
\csep[1.4]{1} \begin{array}{l}
 \Zeta{_1} {x_1} +\Zeta{_2} {x_2}+{x_3} \\
 \Zeta{_1} {x_1} +\Zeta{_2} {x_3}+{x_2} \\
 \Zeta{_1} {x_2} +\Zeta{_2} {x_1}+{x_3} \\
 \Zeta{_1} {x_2} +\Zeta{_2} {x_3}+{x_1} \\
 \Zeta{_1} {x_3} +\Zeta{_2} {x_1}+{x_2} \\
 \Zeta{_1} {x_3} +\Zeta{_2} {x_2}+{x_1}
\end{array}		& 
\begin{array}{ll}
6\,\displaystyle{\prod_{i=1}^{3}}x_{i} &+3\,\Zeta{_1} \left({x_1}{x_3}^2\!+{x_2} {x_1}^2\!+{x_3}{x_2}^2 \right)+\\&+3\, \Zeta{_2} \left({x_1}{x_2}^2\! +{x_2} {x_3}^2\!+{x_3} {x_1}^2\right)+\displaystyle{\sum_{i=1}^{3}}x_{i}^{3} \\\\
6\, \displaystyle{\prod_{i=1}^{3}}x_{i}&+3\, \Zeta{_2} \left({x_1}{x_3}^2\!+{x_2} {x_1}^2\!+{x_3}{x_2}^2 \right)+\\&+3\, \Zeta{_1} \left({x_1}{x_2}^2\! +{x_2} {x_3}^2\!+{x_3} {x_1}^2\right)+\displaystyle{\sum_{i=1}^{3}}x_{i}^{3}
\end{array}	\\\hline
3	& {x_1}+{x_2}+{x_3}	 	&	\left({x_1}+{x_2}+{x_3}\right)^3\\\hline
\end{array}
\end{equation*} 
{The Resolvent equation for the cubic is then of degree 2. }
\begin{equation*}
\scR_n(y)=\left(y-\beta _1\left(x_i\right)\right)\left(y-\beta _2\left(x_i\right)\right)=y^2-b y -\frac{1}{27}a^3
\end{equation*} 
{Its solutions }
\begin{equation*}
\begin{array}{ccc}
 Y_1(a,b)=\frac{b}{2}+\sqrt{\left(\frac{b}{2}\right)^2+\left(\frac{a}{3}\right)^2} & ; & Y_1(a,b)=\frac{b}{2}-\sqrt{\left(\frac{b}{2}\right)^2+\left(\frac{a}{3}\right)^2}
\end{array}
\end{equation*} 
{then imply }
\begin{equation*}
\begin{cases}
 \rho _1=\zeta _1x_1+\zeta _2x_2+x_3 & = 3 Y_1(a,b)^{1/3} \\
 \rho _2=\zeta _2x_1+\zeta _1x_2+x_3 & =3 Y_2(a,b)^{1/3} \\
 \rho _3=x_1+x_2+x_3 & =0
\end{cases}
\end{equation*} 
{with solutions}
\begin{equation*}
\begin{cases}
 x_1&=\zeta _1Y_1(a,b)^{1/3}+\zeta _2Y_2(a,b)^{1/3} \\
 x_2&=\zeta _2Y_1(a,b)^{1/3}+\zeta _1Y_2(a,b)^{1/3} \\
 x_3&=Y_1(a,b)^{1/3}+Y_2(a,b)^{1/3}
\end{cases}
\end{equation*} 
For $n=5$ the smallest number of nontrivial $\beta $ is $24$.
\newline
Lagrange shows first that for any rational expression $\scF\left(x_1,\ldots ,x_n\right)$ the number of its distinct values under
all permutations is $m=\frac{n!}{\left|\scI_{\scF}\right| }$ , where $\left|\scI_{\scF}\right|$ is the number of permutations that leave
$\scF$ invariant. Nowadays Lagrange's theorem states that the order of any subgroup $H$ of a group $G$ divides the order of $G$.
}
\fn{6}{\textbf{\large Aléxandre-Théophile Vandermonde (1735-1796)}
\setlength{\parindent}{0in}
\medskip\newline
 Vandermonde explicitly states that an algebraic expression for the roots of polynomials must be ambiguous given that the enumeration
of roots is arbitrary.
\newline
For a polynomial \begin{equation*}
                  \scP_n(x)=x^n+a_1x^{n-1}+\ldots +a_{n-1}x+a_n
                 \end{equation*}
an expression $\scF\left(x_1,\ldots ,x_n\right)$ involving radicals
$\sqrt[n]{}$ exists such that, depending on the choice of the radical, each root $x_k$ is determined. 
Furthermore this expression $\scF\left(x_1,\ldots ,x_n\right)$ can be written solely in terms of the coefficients $\mathfrak{a}=\left\lbrace a_1,a_2,\ldots ,a_n\right\rbrace $ given that
\begin{equation}
\begin{array}{lcccc}
 a_1 = \displaystyle-\sum_{i}\, x_i\;, &
 a_2 =\displaystyle {\sum_{i<j} }\,x_{i }\,x_j\;, &
 a_3 =\displaystyle -{\sum_{i<j<k} }\,x_i\,x_j\,x_k\;, &
 \dots &
 a_n =\displaystyle (-1)^n\prod_{i}\,  x_i
\end{array}\label{NWT}
\end{equation} 
Finally, $\scF\left(x_1,\ldots ,x_n\right)$ is in fact invariant for any permutation of the roots.
\newline
{For $n\leq 4$ Vandermonde uses explicitly}
\begin{equation*}
\scF\left(x_1,\ldots ,x_n\right)=\frac{1}{n}\sum _{i=1}^n x_i+\sum _{k=1}^{n-1} \beta _k\left(x_1,\ldots ,x_n\right)^{1/n}
\end{equation*} 
{where the $\beta _k\left(x_1,\ldots ,x_n\right)$ are the Lagrange resolvents}
\begin{equation*}
\beta _k\left(x_1,\ldots ,x_n\right)=\left( \frac{1}{n}\sum _i \left(\Zeta{{}_k}^{(n)}\right)^ix_i \right)^n
\end{equation*} 
{For $n=2$, we have $\Zeta{{}_k}^{(2)}=(-1)^k$, and furthermore $x_1+x_2=-a_1$ and $x_1x_2=a_2$ so}
\begin{equation*}
\beta _1\left(x_1,\ldots ,x_n\right)=\frac{1}{4}\left(x_2-x_1\right)^2=\frac{1}{4}\left(a_1^2-4a_2\right)
\end{equation*} 
{We thus obtain}
\begin{equation*}
\scF\left(x_1,x_2\right)=\frac{1}{2}\left(x_1+x_2\right)+\left(\frac{1}{4}\left(x_2-x_1\right)^2\right)^{1/2}=\frac{1}{2}\left(-a_1+ \left(\zeta_1^{(2)}\right)^j\sqrt{a_1^2-4a_2}\right)
\end{equation*} 
In view of the difficulties of his method for $n>4$ Vandermonde studies then the conditions for which higher degree equations may
have solutions with radicals, which led them to study the cyclotomic equations $x^p-1=0$.
\newline
Since for $p$ non-prime its roots may be expressed by radicals if the roots of its prime factors are, all that remains is the study of
cyclotomic equations with $p$ prime. Since $\zeta =1$ is always a root, one must study then the roots of
\begin{equation}
\frac{x^p-1}{x-1}=x^{p-1}+x^{p-2}+\ldots +x+1=0
\label{VD}                                                
\end{equation} 
{In his method Vandermonde proposes that one finds for a cyclotomic equation $x^p-1=0$ a smaller set of $m<p$ symmetric rational
functions $\beta _i\left(\zeta _1,\ldots ,\zeta _p\right)$ of the roots of the original equation such that the $\beta _i$ can be obtained
as roots of algebraic equations of degree $m$.}
\newline
{For $p$ prime, the division of expression (\ref{VD}) by $x^q$, with $q=\frac{p-1}{2}$, and a variable change
to $Y=x+\frac{1}{x}$ expresses equivalent equations of lower degrees.}
\begin{equation*}\renewcommand\arraystretch{1.5}\fontsize{6}{9}
\hskip -12pt\begin{array}[c]{|c|c|}\hline
 1+x+x^2=0 & 1+Y=0\\\hline
 1+x+x^2+x^3+x^4=0 & -1+Y+Y^2=0\\\hline
 1+x+x^2+x^3+x^4+x^5+x^6=0 & -1-2 Y+Y^2+Y^3=0\\\hline
 1+x+x^2+x^3+x^4+x^5+x^6+x^7+x^8=0 &  1-2 Y-3 Y^2+Y^3+Y^4=0 \\\hline
 1+x+x^2+x^3+x^4+x^5+x^6+x^7+x^8+x^9+x^{10}=0 & 1+3 Y-3 Y^2-4 Y^3+Y^4+Y^5=0\\\hline
\end{array}
\end{equation*}
Vandermonde is the first to find solutions by radicals of the cyclotomic equation $x^p-1=0$ for $p=11$, and as a consequence
obtains such solutions for a reduced equation of degree $\frac{p-1}{2}=5$.
}

At the dawn of the XIX$\rm  {}^{th} $ century Paolo Ruffini in his ``Teoria Generale delle Equazioni'' attempts a $ 516 $ page demonstration that general equations of degree at least $ 5 $ are not solvable by radicals. In $ 1810$, in view of a poor reception by the mathematical community, Ruffini submits an improved version to the French academy, but the referees Lagrange, Lacroix and Legendre were so delayed with an answer that Ruffini wrote to the president of the academy to withdraw the submitted work.

In $ 1801 $ Gauss studies cyclotomic equations and achieves important results regarding the solvability of these equations that would serve as foundation for the theory of algebraic equations. In the process he solves a $ 2\ 000 $ year-old problem on the construction of a polygon with straight edge and compass.
\fn{7}{\textbf{\large Carl Friedrich Gauss (1777-1855)}
\setlength{\parindent}{0in}
\medskip\newline
\begin{minipage}{0.2\linewidth}
 \includegraphics[scale=.35,keepaspectratio=true]{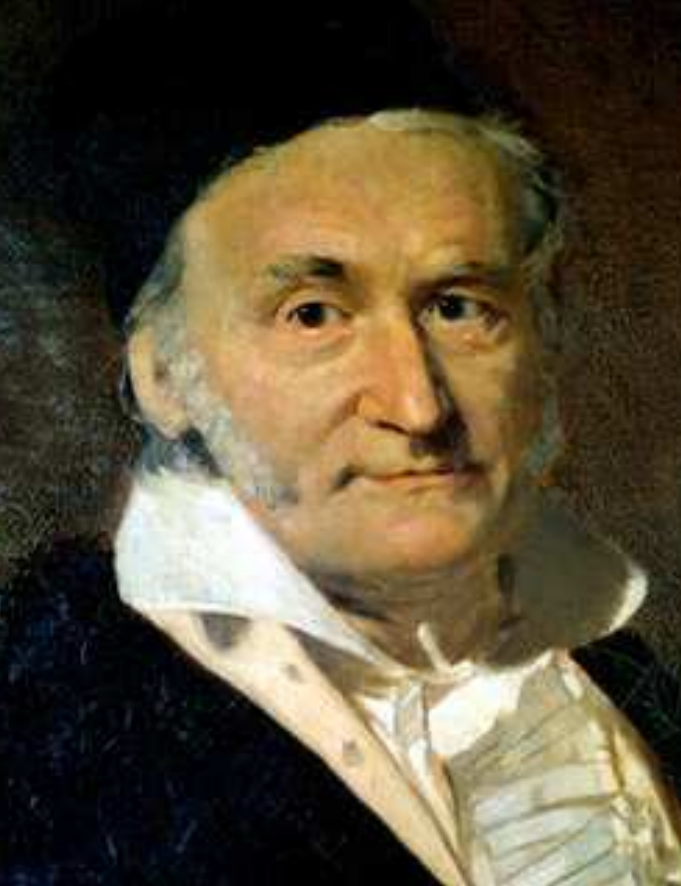}
\end{minipage}\hspace{.05\linewidth}
\begin{minipage}{0.75\linewidth}
In 1801 Gauss studies cyclotomic equations (i.e. equations for the division of the circle)  $ x^p-1= 0 $.
He finds that the roots of $x^p-1=0$ can be rationally expressed in terms of the roots of a sequence of equations $\mathfrak{Z}_m(x)=0$ whose degrees $m$ are prime factors of $p-1$, and whose coefficients are rational in the roots of the preceding equations in the sequence. 
\begin{equation*}
\begin{array}{cc}
x^p-1=\displaystyle\prod _{m|p} \mathfrak{Z}_m(x)           &  \\
\mathfrak{Z}_m(x)=\displaystyle\prod _i^{\varphi (m)} \left(x-\zeta
_i^{(m)}\right)=\displaystyle\prod _{k|m} \left(1-x^{\frac{m}{k}}\right)^{\mu (k)} & 
\end{array}
\end{equation*}
\end{minipage}
\medskip\newline
{where $\varphi (m)$ is the Euler function (counting positive co-primes $i$ of $m$) and $\mu (k)$ is the Möbius function}
\begin{equation*}
\mu (k)=\left\lbrace\begin{array}{cl}
 1 & k=1 \\
 0 & k\, \mbox{has repeated prime factors.} \\
 (-1)^n & k\, \mbox{has all n prime factors distinct.}
\end{array}\right.
\end{equation*} 
{This result shows that some equations of high degree $n$ can be solved by radicals if $n$ is a factor of $p-1$.}
\newline
{The geometric construction problem: if $p-1$ contains no factors other than $2$, a polygon of $p$ sides can be constructed
with straight edge and compass. This is because each of the equations $\mathfrak{Z}_i(x)=0$ is of degree $2$ and each of its roots is
so determined. Thus Gauss showed that all polygons of order $p=2^{2^n}+1$ are constructible, and for the first time in $ 2\ 000 $ years gives the ruler-and-compass construction of a $ 17 $-side polygon. Furthermore he shows that the next higher polygon so constructible would have $ 257 $ sides, then $ 65\ 537 $, etc., as in the sequence}
\begin{equation*}
 2^{2^n}+1=
 3,\, 5,\, 17,\, 257,\, 65\ 537,\, 4\ 294\ 967\ 297,\,\dots
\end{equation*} 
}

Then in 1824 Niels-Henrik Abel (1802-1829) shows that the general algebraic equation of the fifth degree has no solution via radicals. Soon after he extends this result to the nonexistence of a general solution with radicals for algebraic equations of degree greater than four. 

When he dies in $ 1829 $ he was addressing the problem of recognizing if a
particular algebraic equation of high degree can be solvable by radicals but he
couldn't complete his work. Still he was able to show that equations whose roots
$x_k $are rational functions $x_k=F_k\left(x_1\right)$ of a single root $x_1$
all of which verify
$F_k\left(F_j\left(x_1\right)\right)=F_j\left(F_k\left(x_1\right)\right)$ are
solvable. These equations are now called ``Abelian Equations''.  It is Abel that
introduces the concepts of \textbf{Field} and of \textbf{Irreducible
Polynomial}\footnote{A polynomial over a field
$\mathbb{F}$ is said to be \textbf{Reducible} if it can be expressed as a
product of two polynomials of lesser degree over the same field. Otherwise it
will be called \textbf{Irreducible}.} over a field. From Gauss, Ruffini and
Abel's work we conclude that cyclotomic equations of any degree are solvable by
radicals.

Then comes Galois, who submitted  in $ 1829 $, at the young age of $ 18 $, two papers on the solution of algebraic equations to the Academy of Sciences, which were lost by Cauchy. In $ 1830 $ he again presented a paper on his research to the Academy of Sciences. This was sent to Fourier, who unfortunately died soon after and that paper was also lost. The $ 1831 $ article  submitted to Poisson and entitled "Sur les conditions de résolubilité des équations par radicaux" was returned by Poisson as unintelligible.
On the eve of his deadly duel in $ 1832 $ he writes to his friend August Chevalier an account of his researches.

The goal of Galois work was to provide a structure to the search for solutions of algebraic equations that could be expressed by rational functions and radicals involving their coefficients. Around the $ 1830 $'s, and in view of the recent demonstration by Abel ($ 1826 $)  of the nonexistence of such general resolvent formulas for roots of algebraic equations of degrees higher than four, Galois addressed the issue of determining which of these higher order equations are solvable by radicals. The theoretical structure that would decide that this is possible is the \textbf{Galois group of an equation}, and this group must have a structure of what is now known as a \textbf{solvable group}. But in doing so he achieved much more, since his method can be extended to fields other that the rationals and the structure that emerges from the groups can be used to classify field extensions and their subfields.

Galois discovered a method of finding the group of a given equation, the successive partial resolvent equations and their associated groups that result from extending the field of coefficients with the roots of these resolvents. These groups turn out to be subgroups of the original group.
Galois shows that when the group of an equation with respect to a given field is the identity, then the roots of the equation are members of that field.

Application of Galois' theory to the solution of polynomial equations by rational operations and radicals then follows.
When the partial resolvent that serves to reduce to a subgroup $G_2$ the group $G_1$ of an equation is of the form of a binomial equation $x^p=a$ with $p$ prime, then $G_2$ is a \textbf{normal subgroup} of index $ p $ of $G_1$.
Conversely a normal subgroup $G_2$ of prime index $p$ of $G_1$ yields a binomial resolvent $ x^p=a $.

The basic idea is to show that for each algebraic equation of degree $ n $ there is a Galois resolvent $ V_{1}=V_{1}(x_{1},\dots,x_{n}) $ which is a rational expression of its roots $ x_{i} $ with the property that we can in principle find rational expressions for the $ x_{i} $ as 
\begin{equation}
x_{i}=f_{i}(V_{1})
\label{GR}
\end{equation} 
The Galois Resolvent is found to be one of the roots $ V_{k} $ of a polynomial equation $ \scR(V)=0 $, called the \textbf{minimal polynomial}, such that substitution of any of its roots in (\ref{GR}) yields a permutation of the roots
\begin{equation*}
 x_{i_{k}}=f_{i}(V_{k})
\end{equation*} 
These permutations form the Galois Group of the equation.

A basic example illustrates the method. Consider the quadratic equation in $ x $
\begin{equation*}
 x^{2}+a\,x+b=0
\end{equation*} 
with two roots $ x_{1} $ and $ x_{2} $. The Galois resolvent in this case is simply 
\begin{equation*}
V=x_{1}-x_{2}
\end{equation*} 
Taking into consideration the general relations between the roots
\begin{equation}
\displaystyle\sum_{i=1}^{2}x_{i}=x_{1}+x_{2}= -a\qquad;\qquad
\displaystyle\sum_{i<j}^{}x_{i}x_{j}=x_{1}x_{2}= b
\label{NR}
\end{equation} 
we can write $ x_{2}=-x_{1}-a $ and
\begin{equation*}\renewcommand{\arraystretch}{2}
\left\lbrace \begin{array}{l}
\displaystyle x_{1}=\dfrac{V-a}{2}=f_{1}(V)\\
\displaystyle x_{2}=-\dfrac{V+a}{2}=f_{2}(V)
\end{array}\right.
\end{equation*} 
Now from the second relation (\ref{NR}) we get 
\begin{equation*}
 x_{1}x_{2}= b \qquad \Rightarrow\qquad f_{1}(V)\,f_{2}(V)=-\dfrac{V^{2}-a^{2}}{4}=b
\end{equation*} 
and the minimal polynomial is
\begin{equation*}
 \scR(V)=V^{2}+4\,b-a^{2}=(V-\sqrt{a^{2}-4\,b})(V+\sqrt{a^{2}-4\,b})
\end{equation*} 
therefore its roots are $ V_{1}=\sqrt{a^{2}-4\,b} $ and $ V_{2}=-\sqrt{a^{2}-4\,b} $. If we choose to set
\begin{equation*}\renewcommand{\arraystretch}{2}
\left\lbrace \begin{array}{l}
\displaystyle x_{1}=f_{1}(V_{1})=\frac{1}{2}\left(\sqrt{a^{2}-4\,b}-a\right)\\
\displaystyle x_{2}=f_{2}(V_{1})=-\frac{1}{2}\left(\sqrt{a^{2}-4\,b}+a\right)
\end{array}\right.
\end{equation*} 
then it is evident that
\begin{equation*}\renewcommand{\arraystretch}{2}
\left\lbrace \begin{array}{l}
\displaystyle x_{1}^{\prime}=f_{1}(V_{2})=x_{2}\\
\displaystyle x_{2}^{\prime}=f_{2}(V_{2})=x_{1}
\end{array}\right.
\end{equation*} 
Thus the Galois group of the quadratic polynomial is simply a two element set
\begin{equation*}
 G=S_{2}=\left\lbrace \sigma_{1}=(1)(2)\quad,\quad\sigma_{2}=(12)\right\rbrace
\end{equation*} 
The only normal subgroup of this group is trivially the identity permutation $ \sigma_{1} $, whose index in $ G $ is $ 2 $, a prime number. 

That this procedure can theoretically be implemented for any algebraic equation, and its group thus defined, is the genius of Galois theory, but to carry this out explicitly is mostly highly impractical, as Galois himself already knew. But if there are other ways to get to the group of the equation, then its solvability by radicals can then be ascertained solely by looking at the subgroup structure of this group, and that is oftentimes quite easy to do. In this category are the methods involving the coefficient field extension which allows for a definition of partial resolvents in a new field, and the construction of the composition series as a result. The example in Note 8 below illustrates this point. 
\fn{8}{\textbf{\large Évariste Galois (1811-1832)}
\setlength{\parindent}{0in}
\medskip\newline
\begin{minipage}{0.2\linewidth}
 \includegraphics[scale=.3,keepaspectratio=true]{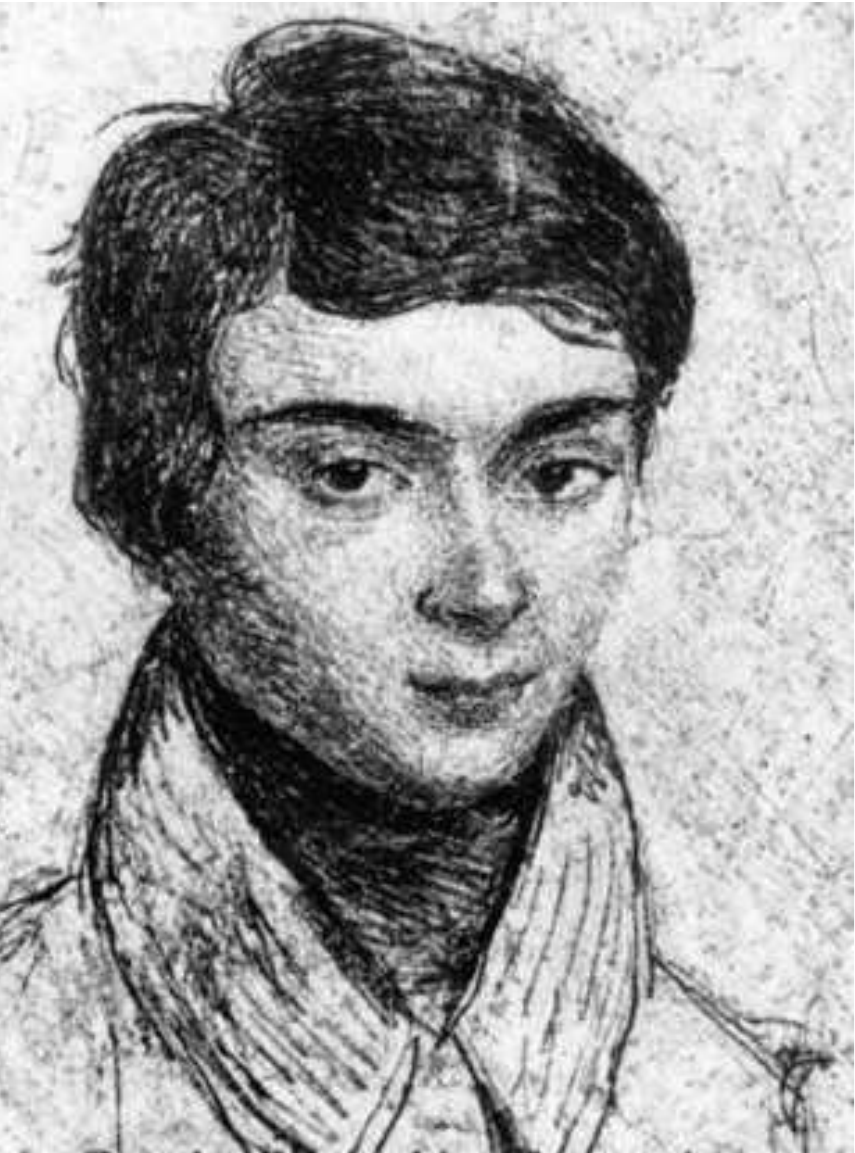}
\end{minipage}\hspace{.05\linewidth}
\begin{minipage}{0.75\linewidth}
 In the following example we illustrate from a known solvable equation the concepts developed by Galois. Of course the theory works without knowing first the solutions of the equations, and it makes use of theorems relating symmetric functions of the roots, the construction of the Galois partial resolvents and its use in determining a composition series of subgroups. We can also see how the concept of field extension enters the theory. Notice that even though the theory shows the existence of Galois Resolvents and how it relates to the construction of the Galois Group of the equation, it is most of the times extremely hard to find these functions explicitly, therefore an indirect approach has to be taken. Galois himself knew that his method was not a way to construct explicit representations of the roots of polynomial equations via rational expressions of the coefficients in a particular extended field (what nowadays would be called the splitting field).
\end{minipage}
\medskip\newline
Consider the quartic (biquadratic) equation
\begin{equation}
\scP_4(x)=x^4+a x^2+b=0   
\label{EQ1}                    
\end{equation} 
In this case we already know the roots of this equation since we can solve it for $ y=x^{2} $ as a quadratic equation and then use $ x_{i_{\pm}}=\pm\,\sqrt{y_{i}} $. One possible labeling of these solutions would be 
\begin{equation*}
  \begin{array}{ccc}
x_{1}=-x_{2}=\sqrt{y_{1}}\,&;& x_{3}=-x_{4}=\sqrt{y_{2}}                                                               
\end{array}
\end{equation*} 
but there are $ 23 $ other choices.
We will find out the Galois group of equation (\ref{EQ1}) from certain symmetries associated with the arbitrary choice of indexes for the roots.
\newline
Set $\mathbb{Q}_1(a,b)$ the field of rational expressions in $a,b$ with coefficients in the field $ \mathbb{Q} $ of rational numbers. The field $\mathbb{Q}_1(a,b)$ is called an \textbf{extension} of the  rational field $ \mathbb{Q} $. Then the following relations hold in $\mathbb{Q}_1(a,b)$
\begin{equation}
\begin{array}{ccc}
 x_1+x_2=0 & ; & x_3+x_4=0
\end{array}
\label{Q1}
\end{equation} 
Now the root set of equation (\ref{EQ1}) is invariant under the $4!=24$ permutations of their labeling, but only $8$ of these permutations will leave the relations (\ref{Q1}) invariant, and they are in the set $ G_{1} $:
\begin{equation*}
G_1=\left\lbrace
\begin{array}[c]{llll}
 \sigma _1=(1)(2)(3)(4), & \sigma _2=(12)(3)(4), & \sigma _3=(1)(2)(34), & \sigma _4=(12)(34) \\\\
\sigma _5=(13)(24), & \sigma _6=(1423), & \sigma _7=(1324), & \sigma _8=(14)(23) 
\end{array}
\right\rbrace
\end{equation*} 
This set is a group under the composition of permutations $ \sigma_{i} $.
\footnote{ In fact $ G_{1} $ is a \textbf{normal subgroup} of $S_4$ of \textbf{prime index} $ 3 $ (because $ \frac{|S_{4}|}{|G_{1}|}=\frac{24}{8}=3 $).} 
This group $ G_{1} $ is the \textbf{Galois group} of the equation (\ref{EQ1}) as it is the largest subgroup of the symmetric group $ S_{4} $ that leaves invariant the basic set (\ref{Q1}) of rational functions of the roots with coefficients in $ \mathbb{Q}_{1}(a,b) $. 
\newline
Now we know from Newton's relations (\ref{NWT}) that $ {x_{1}}^{2}+{x_{3}}^{2}=-a $ and $ {x_{1}}^{2}{x_{3}}^{2}=b $ so the next relation 
\begin{equation}
\begin{array}{ccccc}
{x_1}^2-{x_3}^2=\xi_{1}=\sqrt{a^2-4b}
\end{array}
\label{Q2}
\end{equation} 
is not rational in the field $ \mathbb{Q}_{1}(a,b) $, but it is in the field $\mathbb{Q}_2\left(a,b,  \xi_{1}\right)$ which is by definition the extension of $\mathbb{Q}_1(a,b)$ to a field of rational expressions in $a,b$ and $ \xi_{1}$.
\newline
It should be apparent that now only the first $ 4 $ permutations in $G_1$ leave these relations (\ref{Q1}) and (\ref{Q2}) invariant:
\begin{equation*}
G_2=\left\lbrace\, \begin{array}[c]{llll}
\sigma _1=(1)(2)(3)(4),& \sigma _2=(12)(3)(4), & \sigma _3=(1)(2)(34), & \sigma _4=(12)(34)
\end{array} \right\rbrace 
\end{equation*} 
since ${x_1}^2={x_{2}}^2$ and ${ x_{3}}^2={x_{4}}^2 $.
\newline
Now $ G_{2} $ is a \textbf{normal} subgroup of $ G_{1} $ with prime index $ 2 $. Notice also that the $ \xi_{1} $ in equation (\ref{Q2}) is in fact a solution of the \textbf{partial resolvent} polynomial of degree $ \frac{|G_{1}|}{|G_{2}|}=2 $ in $ \mathbb{Q}_{1}(a,b) $
\begin{equation}
\xi^{2}-a^2+4b=0
\label{PR1}
\end{equation} 
From the relations (\ref{Q1}) we can also derive a new expression
\begin{equation}
x_1-x_2=2\,\xi_{2}
\label{Q3}
\end{equation} 
which we can view\footnote{Squaring both sides of (\ref{Q3}), adding ${x_{3}}^{2} $ and noting that $ {x_1}^{2}={x_2}^{2}={\xi_{2}}^{2} $ and $ \xi_1={x_1}^2-{x_3}^2 $, we get ${x_1}^{2}-2\,x_1\,x_2+{x_{2}}^{2}+{x_{3}}^{2}=-a+2\,{\xi_{2}}^{2}+{x_{1}}^{2}=4\,{\xi_{2}}^{2}+{x_{3}}^{2} $ and thus (\ref{PR2})} as a root of the partial resolvent polynomial of degree $ 2 $ in $ \mathbb{Q}_{2}(a,b,\xi_{1}) $
\begin{equation}
 2\,\xi^{2}+a-\xi_{1}=0
\label{PR2}
\end{equation} 
and thus the relation (\ref{Q3}) is rational in the extended field $\mathbb{Q}_3\left(a,b, \xi_{1},\xi_{2}\right)$. 
Now the group leaving invariant all previous root relations plus this one (\ref{Q3}) is
\begin{equation*}
G_3=\left\lbrace
\begin{array}[c]{ll}
 \sigma _1=(1)(2)(3)(4), & \sigma _2=(12)(3)(4) 
\end{array}
\right\rbrace
\end{equation*} 
This group $ G_{3} $ is also \textbf{normal} in $ G_{2} $ with prime index $ 2 $, as expected from the degree of the partial resolvent equation (\ref{PR2}).
\newline
Likewise, the relation
\begin{equation}
x_3-x_4=2\,\xi_{3}
\label{Q4}
\end{equation}
is a root of the partial resolvent polynomial equation of degree $ 2 $
\begin{equation}
 2\,\xi^{2}+a+\xi_{1}=0
\label{PR3}
\end{equation} 
and then only the identity permutation leaves all these root relations (\ref{Q1}), (\ref{Q2}), (\ref{Q3}) and (\ref{Q4}) invariant in the extension field $\mathbb{Q}_4\left(a,b, \xi_{1},\xi_{2},\xi_{3}\right)$. Thus
\begin{equation*}
G_4=\left\lbrace
\begin{array}[c]{l}
  \sigma _1=(1)(2)(3)(4)
\end{array}
\right\rbrace\equiv\mathit{1}
\end{equation*} 
is trivially a normal subgroup of $ G_{3} $ with prime index $ 2 $ too. 
\newline
We have thus obtained a sequence of subgroups
\begin{equation*}
\mathit{1}\vartriangleleft G_{3}\vartriangleleft G_{2}\vartriangleleft G_{1}
\end{equation*} 
called a \textbf{composition series} of the Galois group of prime indexes $ 2:2:2 $. The fact that all the above resolvent equations are \textbf{binomial equations} $ x^{p}-A=0 $ with $ p $ prime is intimately connected with the fact that the equation is solvable by radicals. The fact that all the indexes in the composition series of a particular group are prime numbers establishes the group as \textbf{solvable}.
}

If all the partial resolvents are binomial equations then, since Gauss has shown that binomial equations can be solved by radicals, it follows that the original equation can be solved by radicals by extending the original field of the equation coefficients to one with all the roots added by successive adjunction of radicals.
Conversely an equation solvable by radicals must have partial resolvent equations which are binomial equations of prime degree.

The theory of solvability by radicals is thus equivalent to demanding that the successive subgroups of the original Galois group (the \textbf{composition series}) must each be a \textbf{maximum normal subgroup} of the preceding group.
Thus, in order for an equation to be solvable by radicals the indexes in the composition series must be a sequence of prime numbers, and a group $G$ which contains a composition series of prime indexes is said to be \textbf{solvable}.

For instance, the group $S_4$ is solvable, having a composition sequence 
\begin{equation*}
\mathit{1}\vartriangleleft C_{2}\vartriangleleft V_{4} \vartriangleleft  A_{4}\vartriangleleft S_{4}
\end{equation*}  where $ A_{4} $ is the alternating normal subgroup ($ 12 $ even permutations), 
\begin{equation*}
 V_{4}=\left\lbrace (1)(2)(3)(4),\,(12)(34),\,(13)(24),\,(14)(23)\right\rbrace
\end{equation*}
is the Klein 4-group or Vierergruppe, 
\begin{equation*}
C_{2}=\left\lbrace (1)(2)(3)(4),\,(12)(34)\right\rbrace
\end{equation*} is the Cyclic 2-group and $ \mathit{1} $ the identity group, with indexes $2:3:2:2$. Thus, since $ S_{4} $ is the Galois group of the general algebraic equation of fourth order, these equations are solvable. 

However, for general algebraic equations of order $ n>4 $, the Galois groups $ S_{n} $ are never solvable, since their maximal normal subgroup is the alternating group $ A_{n} $ of order $ \frac{n!}{2} $, and this type of group for $ n>4 $ does not have any normal subgroups itself apart from the identity. Thus the composition series always has indexes $ 2:\frac{n!}{2} $, and $ \frac{n!}{2}  $ is never prime for $ n\geq 4 $.

Only in $ 1846 $ will Liouville, reading through Galois' work in possession of the brother Alfred Galois, finally understand the magnitude of his findings and publish them in the Journal de Mathématiques. For the next 15 years Cayley, Dedekind and Kronecker will work on the subject. In $ 1854 $ Cayley generalizes the concept of permutation group of a finite set and defines a \textbf{Finite Abstract Group} as any finite set with an associative composition and a neutral element. His work is largely ignored by the community at the time as matrices and quaternions were new and not well known. 

By $ 1866 $ Serret lectures at the Sorbonne about Galois's work after Liouville publication in the Journal de Mathématiques.
It is said that Serret's 3$\rm {}^{rd} $ edition of the Cours d'Algébre Supérieure was so popular that its adoption in France's mathematics curricula for the next $ 50 $ years hindered there the divulgation of the latest developments in group theory.
Serret first studies representations of substitutions by transformations of the form
\begin{equation*}
  z\to z'=\frac{p\, z+q}{p'z+q'}\quad \text{with}\quad p\,q'-p'q\neq 0
\end{equation*} 
Then in $ 1868 $ Jordan starts the first investigations of Infinite Groups with paper ``Mémoire sur les groupes de mouvements'' after Bravais in $ 1849 $ had studied groups of motions to determine the possible structure of crystals. 

In $ 1870 $ Jordan publishes the "Traité des substitutions et des équations algébriques" where he organizes his work on
finite substitution groups and their connection with Galois theory.
He is the first mathematician to focus on the group theory aspects of the work rather than on finding the roots of equations.
Jordan solves Abel's problem of determining which equations of a given degree are solvable by radicals.
He concludes that the groups of such solvable equations are commutative, which he then names "Abelian".

In $ 1873 $ Jordan introduces the notion of \textbf{quotient group}. He establishes that on different composition series of the same group there are always the same number of elements and the order of the quotient groups is the same up to ordering.
In $ 1889 $ Hölder shows that the quotient groups in a composition series are isomorphic up to ordering (Jordan-Hölder Theorem).

It is Jordan who initiates the study of geometric transformations with groups as he studies infinite groups of translations and rotations and represents substitutions by linear transformations of the form
\begin{equation*}
x_{i }\to x_i^{ \prime }=\sum _{j=1}^n A_{ij}x_j
\end{equation*} 

In the late XIX$\rm  {}^{th} $ century, German mathematics was very active pursuing problems in Number Theory. Richard Dedekind ($ 1831-1916 $)
developed the foundations of modern Galois Theory through its applications in field theory. 

Leopold Kronecker ($ 1823-1891 $) 
used group theory not in terms of permutations of roots of an equation but as a group of automorphisms of the coefficient
field and its extensions.

In $ 1893 $ Heinrich Weber ($ 1842-1913 $)
presents Galois Theory in terms of group and field theory, using theorems no longer restricted to rationals but applied to arbitrary fields, and in $ 1895 $ in his ``Lehrbuch der Algebra'' explicitly extends the notion to infinite groups. The chapter on Galois Theory in Weber's algebra textbook makes no reference to the solvability of groups or the solvability of equations by radicals. From this point on Galois Theory was no longer concerned with the practical construction of roots of algebraic equations.
Certain problems in Field theory can now be reduced to group theory through Galois Theory.

The modern-day Galois Theory is understood in the formulation of Emil Artin ($ 1898-1962 $).
In $ 1930 $ Artin finally established Galois Theory as the study of Field Extensions and their Automorphism Groups, revealing the parallelism in their structure.

Artin abandoned the approach of building a sequence of field extensions by adjoining resolvents to the coefficient field and introduces instead the \textbf{Splitting Field} of the equation as the smallest field containing the roots and the coefficients.

In $ 1931 $ Bartel van der Waerden publishes "Moderne Algebra" using the lectures by E. Artin and E. Noether. 
In $ 1963 $ the Feit-Thompson theorem is proven showing that a finite group of odd order is necessarily solvable.
A very interesting account of the development of Galois theory up to Artin's work can be found in \cite{BMK71}.

In the meantime, since the mid-$ 1800 $ the group concept irrespective of its Galois connotations found its way to applications in Differential Geometry through the works of Sophus Lie, Felix Klein and many others. Klein's Erlangen Program is the best statement of the universality of the group notion, and its use via Representation Theory in Physics is the trademark of the XX$\rm  {}^{th} $ century physics, were we associate physical symmetry groups and their irreducible representations with measurable characteristics of natural systems.
Using either discrete or continuous, differentiable groups and its infinitesimal counterpart, the Lie Algebra, one can understand (or model) such diverse aspects as Spontaneous Symmetry Breaking and Renormalization Methods \cite{JPE84},
the reason for the structure of the Periodic Table \cite{SS95},
the Selection Rules in Quantum Mechanics \cite{JPE84} \cite{SS95} \cite{MH89},                                   
the Relativistic and Non-relativistic Dynamical Symmetries \cite{AOB86},
the Symmetries of Differential Equations and Conservation Laws \cite{DHS93} \cite{HS90} \cite{PJO00},
the reason why Parity is not a natural symmetry \cite{SS95}
or why all Relativistic Wave Equations are an expression of projection operators of Induced Representations of the Lorentz or Poincaré groups \cite{JPE84}\cite{AOB86}.

Coming back to Mira Fernandes,  we can see how close he was in his $ 1910 $ thesis to the formal developments that the finite abstract group theory had brought in the last decade of the XIX$\rm  {}^{th} $ century. The example groups that he studies (Abelian, Metacyclic, General linear, Modular and Symmetric) are not connected in this work to any applications to the resolution of algebraic equations. The results that he exposes are of general nature in finite group theory, but the name ``Galois Theories-I'' in his thesis indicate already his intention of applying these group theoretical concepts to the theory of algebraic equations, a work that he finishes with his $ 1931 $ initiation text on algebraic solvability of equations in the book entitled ``Grupos de Substituições e Resolubilidade Algébrica''\footnote{Substitution Groups and Algebraic Solvability}\cite{AdMF29b}. A couple of years before, in $ 1929 $, he had reedited most of his thesis material as part I of this book \cite{AdMF29a}. 

Mira Fernandes' approach to Galois groups and algebraic equations in his $ 1931 $ work \cite{AdMF29b} is a little different in that he introduces the fundamental function 
\begin{equation*}
 \phi_{1}=\sum_{i=1}^{n}\alpha^{i-1}x_{i}
\end{equation*} with $ \alpha $ an arbitrary constant not a root of the discriminant of his ``Galois Resolvent'' 
\begin{equation}
 R(y)=\prod_{i=1}^{n!}(y-\phi_{i})=0
\label{GR1}
\end{equation} 
where the $ \phi_{i} $ are $ n! $ expressions obtained from permutations of the roots $ x_{i} $ in the fundamental function $ \phi_{1} $. He then defines the Galois Group of the algebraic equation $ f(x)=0 $ as the set of permutations  of the roots $ x_{i} $ that transform $ \phi_{1} $ in each of the distinct roots $ \phi_i $ of an irreducible component of (\ref{GR1}).
He starts by showing that such a resolvent always exists, and proceeds to the theoretical construction of the resolvent polynomial and to prove the theorems in group theory that are now associated with the interpretation of solvability by radicals of algebraic equations. He treats explicitly the general cases for Polynomial equations, and gives a few examples of Abelian equations.

It is reasonable to assume that by then his interest in this Galois theory is merely academic, and probably motivated by his interest in the life of Évariste Galois, which he addresses in a communication published in 1933 by the centennial commemoration of Galois' death\cite{AdMF33}\cite{AdMF33-a}. His $ 1931 $ work on Galois' theory is then entirely dedicated to the solvability theory of algebraic equations, and little focus is given to the field theoretical applications providing links between groups and field extensions that by then would become the standard when referring to Galois' theory\cite{EA38}. In $ 1936 $ he would once again refer Galois work in the context of the importance of group theory in the communication on the Evolution of Variational Calculus presented on the commemoration of the bi-centennial of Lagrange's birth \cite{AdMF}.

After his thesis work on Galois Theories Mira Fernandes seems to be dedicated mostly to academic activities, but surely became involved in studying the applications of Lie groups and Differential Geometry in the context of connection models for general-relativistic unitary theories (see \cite{Lem09} in this publication), where his production was extensive and internationally recognized after a gap of about $ 15 $ years without publications. 

His early awareness of the breath of implications that group theory concepts brought to the understanding of geometrical and physical theories, as can be glanced by his historical references to Klein's Erlangen Program in his talks, was certainly initiated by his studies in Galois Theories, and because of it he can be cast as one of the early visionaries of the new methodology that would become essential in most areas of XX$\rm  {}^{th} $ century mathematical physics.

\theendnotes
\bibliographystyle{alpha} 
\bibliography{MF}

\end{document}